\documentclass[final,onefignum,onetabnum]{siamart171218}

\usepackage{lipsum}
\usepackage{amsfonts}
\usepackage{graphicx}
\usepackage{epstopdf}
\usepackage{algorithmic}
\usepackage{mathrsfs}
\ifpdf
  \DeclareGraphicsExtensions{.eps,.pdf,.png,.jpg}
\else
  \DeclareGraphicsExtensions{.eps}
\fi

\newcommand{\bs}[1]{\boldsymbol{#1}}

\newsiamremark{remark}{Remark}
\newsiamremark{hypothesis}{Hypothesis}
\crefname{hypothesis}{Hypothesis}{Hypotheses}
\newsiamthm{claim}{Claim}

\title{A note on generating Voronoi cells with a\\ given size distribution\thanks{The authors appreciate support by the US Office of Naval Research (ONR) grant \#N00014-19-1-242. Please contact the authors for comments and questions.}}
\author{Georg Stadler\thanks{Courant Institute of Mathematical
    Sciences, New York University 
  (\email{stadler@cims.nyu.edu}) (\email{gg2924@nyu.edu})}
\and Gonzalo G.~De Diego\footnotemark[2]}

\usepackage{amsopn}

\usepackage{color}              
\usepackage{tikz}
\usetikzlibrary{positioning}
\usepackage{pgfplots}
\usepackage{pgfplotstable}

\usepgfplotslibrary{fillbetween}
\usepgfplotslibrary{colorbrewer}
\usepgfplotslibrary{groupplots}
\usepgfplotslibrary{statistics}
  \pgfplotsset{compat=newest}

\begin{document}

\maketitle

\begin{abstract}
This note describes a simple method to draw random points such that the cells of the corresponding Voronoi tesselation (approximately) satisfy a desired size distribution, for instance, follow a power law. The method is illustrated and numerically verified in two dimensions, and we also provide a simple implementation.
\end{abstract}



\section{Introduction}
The size distribution of Voronoi cells generated from random points drawn from a uniform distribution has been examined both theoretically and numerically \cite{ferenc2007size}. In one dimension, the distribution can be calculated analytically, but analytic expressions are unknown in two or three dimensions. However, approximations that closely match the numerically determined size distributions for large point counts ($10^6$ to $10^7$) are available. In two dimensions, the distribution of cells of area $y$ is effectively modeled by the following two-parameter generalized Gamma function:
\begin{equation}\label{eq:distribution}
f(y) = \frac{b^a}{\Gamma(a)} y^{a-1}\exp(-by).
\end{equation}
Typically used values are $a=3.61$ and $b=3.57$. A similar result holds in three dimensions, but for the rest of note we focus on the two-dimensional case.
There exists a slightly more complicated three-parameter formula that yields better fits for the tails of the distribution \cite{ferenc2007size}.
The probability density function \eqref{eq:distribution} assumes a large number of $N$ uniformly distributed points in a two-dimensional domain of area $N$, so the expectation of this distribution is $\mathcal E[f]=1$.

\begin{figure}[th]
    \centering
    \includegraphics[width=0.5\linewidth]{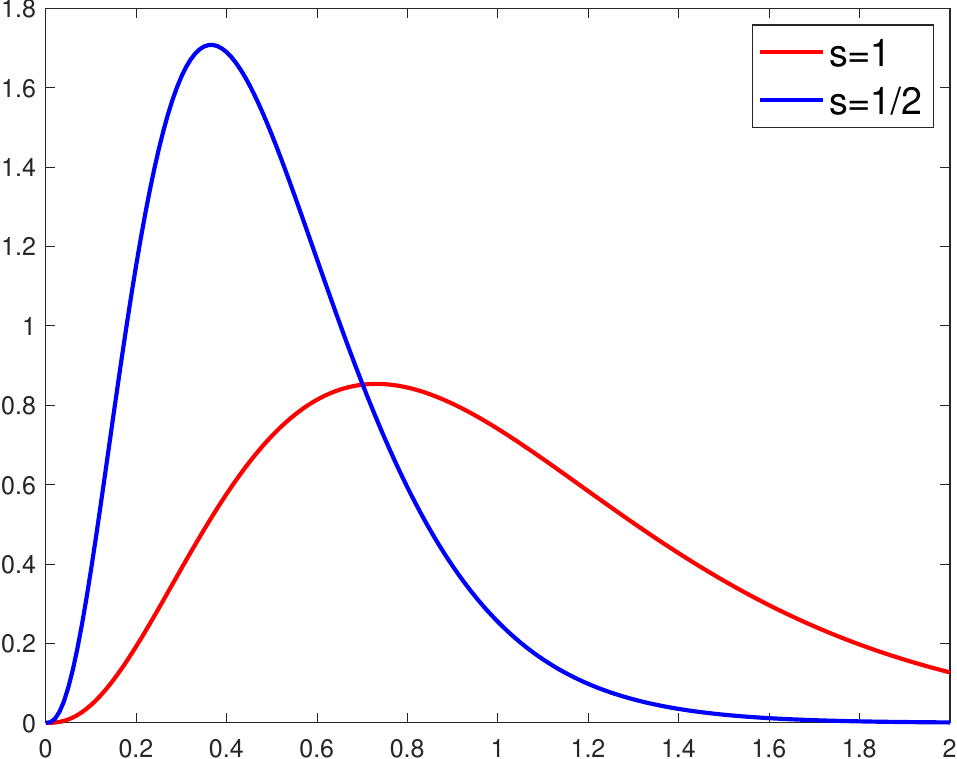}
    \caption{Densities or scaling factors $s=1$ (blue) and $s=1/2$ (red).}
    \label{fig:1}
\end{figure}

Next, we derive how this distribution changes for denser point clouds.
For that purpose, we introduce a scale parameter $s > 0$, which corresponds to the average cell area. That is, on a domain of size $N$ we randomly draw $N/s$ points. Thus, $1/s$ is the average number of points per unit area. 
%
The $s$-dependent distribution $f_s(\cdot)$ then becomes
\begin{equation}\label{eq:s-distribution}
f_s(y) = \frac{b^a}{s\Gamma(a)} \left(\frac{y}{s}\right)^{a-1}\exp\left(-\frac{by}{s}\right).
\end{equation}
The densities for $s=1,1/2$ are shown in Figure~\ref{fig:1}. Note that the expectation of the area for a parameter $s$ is $\mathcal E[f_s]=s$.

\section{Targeting a cell area distribution}
The Voronoi area distribution \eqref{eq:s-distribution} is obtained for $1/s$ randomly drawn points per unit area.  
In some applications, one may want to achieve a desired cell area distribution. For instance, one may be interested in achieving a power-law cell area distribution over a range of areas. To give an example where this is of interest, consider the discrete element model (DEM) code {\it Subzero} \cite{manucharyan2022subzero, Montemuro23}, which uses Voronoi cells to initialize the shape of ice floes. Satellite observation data of sea ice indicates that the floe area typically follows a power-law distribution \cite{buckley2024seasonal,DentonTimmermans22}. Thus, it may be of interest to initialize DEM simulation with ice floes (i.e., Voronoi cells) that satisfy the area distribution observed in nature. To be specific, the power-law distribution $f^{\text{tar}}(\cdot)$ typical for sea ice is
\begin{equation}\label{eq:f^tar}
    f^{\text{tar}}(y) = Cy^{-m},
\end{equation}
where $y$ denotes the flow area, $C>0$ is a constant, and $m>0$ is the exponent. For sea ice, $m$ is typically between $1.75$ and $2$. 
We aim to approximately achieve such a target distribution by
combining Voronoi cells with different point densities.

To do so we fix a vector of allowed scales $\bs s=(s_1,\ldots,s_K)$, $K\in \mathbb N$ and denote by $\bs\alpha:=(\alpha_1,\ldots,\alpha_K)$ the fraction of the domain $\Omega$ covered with cells originating from the point densities corresponding to these scales. Since $\alpha_i$ are fractions, we must ensure that $\alpha_i\ge 0$ for all $i$, and that $\alpha_1+\ldots+\alpha_K=1$. 
Then, for fixed scales $\bs s$, the distribution of Voronoi cells corresponding to $\bs\alpha$ is
\begin{equation}\label{eq:f_alpha}
    f(\bs \alpha;y):=\sum_{i=1}^K \frac{\alpha_i}{s_i} f_{s_i}(y),
\end{equation}
where the factor $1/s_i$ is due to the fact that when cells are smaller on average, more of them fit in a given area.

To target a power-law distribution such as \eqref{eq:f^tar}, it is typically advantageous to use logarithms. For \eqref{eq:f^tar}, we obtain: 
\begin{equation}\label{eq:log(f^tar)}
    \log (f^{\text{tar}}(y)) = A - m \log(y), \quad \text{ with } A:=\log(C)>0. 
\end{equation}
Combining this with the logarithm of \eqref{eq:f_alpha} in a least-squares objective, we obtain the minimization problem:
\begin{equation}
    J(\bs\alpha,A):= \frac12 \int |\log(f(\bs\alpha;y)) - (A-m\log(y))|^2 \,dy.
\end{equation}
Here, in addition to the weight $\bs \alpha$, we treat the unknown constant $A$ as a variable.
Finally, discretizing the area sizes $y$ in a vector of target sizes $\bs y=(y_1,\ldots,y_L)$ results in the following nonlinear least squares problem:
\begin{equation}\label{eq:least-squares-opt}
    \min_
    {\begin{subarray}{c} A\ge 0,\bs \alpha\ge 0\\ \alpha_1+\ldots+\alpha_K=1\end{subarray}} \frac 12 \sum_{l=1}^L \omega_l\left(\log(f(\bs\alpha;y_l)) - (A-m\log(y_l))\right)^2,
\end{equation}
where $\omega_l>0$ are the quadrature weights. This problem can be solved using numerical optimization software.

\section{Examples for solving \eqref{eq:least-squares-opt}}
We use two examples to illustrate the approach, one with $K=4$ and one with $K=5$. 
\paragraph{Example 1} Using the scale and target area vectors
\begin{equation*}
    \bs s =  (1,1/2,1/4,1/8),\quad  \bs y = (0.1,0.15,0.2,0.3,0.5,0.75,1,1.5,2),
\end{equation*}
we aim to achieve an (approximate) power-law cell area distribution with $m=1.75$.
The vector $\bs y$ is chosen non-uniformly to place more emphasis on smaller cell sizes, where we expect a larger number of cells. For the same reason we choose $\omega_l=1$ for all $l$ in \eqref{eq:least-squares-opt}, i.e., we use a collocation method for the target areas $\bs y$.  
Note that the entries in $\bs s$ are the expectations of the distributions \eqref{eq:s-distribution}, which are combined in \eqref{eq:f_alpha}. Thus, the target sizes specified in $\bs y$ cannot be much larger or smaller than the largest or smallest values in $\bs s$. In particular, the smallest expectation is $1/8=0.125$ and thus one cannot hope to achieve a power-law distribution for cells of much smaller area. Hence, the smallest entry in $\bs y$ is chosen as $0.1$. Similarly, the largest entry in $s$ is $1$, and thus we cannot hope to fit a target that focuses on cells much larger than area $1$ (we choose $2$ as the largest value in $\bs y$).

To solve this overdetermined least squares problem \eqref{eq:least-squares-opt}, we use a nonlinear optimization solver that relies on finite differences for gradients (in particular, we use the function \texttt{fmincon} in MATLAB).\footnote{The simple MATLAB implementation can be found in \url{https://github.com/georgst/voronoi-area}.} This optimization yields the area fraction vector
\begin{equation}\label{eq:alpha_opt-ex1}
    \bs\alpha\approx(0.48,0.10,0.145,0.275).
\end{equation}
That is, $48\%$ of the domain should be covered with Voronoi cells on scale $1$, $10\%$ with cells of scale $1/2$, etc. In Section~\ref{sec:vor} we discuss ways to realize a cloud of random points corresponding to \eqref{eq:alpha_opt-ex1}.

\paragraph{Example 2} Next, we target a power-law distribution with $m=1.5$ and use
\begin{equation*}
    \bs s =  (1,1/2,1/4,1/8,1/16),\quad     \bs y = (0.05,0.075,0.1,0.15,0.2,0.3,0.5,0.75,1,1.5,2).
\end{equation*}
Solving \eqref{eq:least-squares-opt} as described in the previous example, we obtain
\begin{equation}\label{eq:alpha_opt-ex2}
    \bs\alpha\approx(0.50,0.04,0.19,0.11,0.16).
\end{equation}
If instead we use the vector $\bs y$ from Example 1, i.e., do not include the cell areas $0.05$ and $0.075$ in the target vector, we obtain
\begin{equation*}
    \bs\alpha\approx(0.56,0.07,0.15,0.22,0),
\end{equation*}
i.e., compared to \eqref{eq:alpha_opt-ex2}, no cells from the distribution with $s=1/16$ are included.

\section{Choosing Voronoi points for given $\bs \alpha$}\label{sec:vor}
Given an area fraction vector $\bs \alpha$, the simplest way to obtain a distribution of Voronoi cells in a domain $\Omega$ is to split the domain into disjoint subdomains $\Omega_i$ such that the areas satisfy $|\Omega_i|/|\Omega|=\alpha_i$ for all $i$. In each $\Omega_i$, points are then randomly drawn from a uniform distribution such that the number of points per unit of area is $1/s_i$. Although the resulting Voronoi cell areas will (approximately) satisfy the targeted size distribution, one is typically interested in sufficient mixing of the differently sized Voronoi cells.

Thus, one may aim to find a point cloud that is locally denser according to the area fraction vector $\bs \alpha$. However, one would like to avoid, to the extent possible, regions of dense points next to much less dense regions. Such internal boundaries at which the point density changes can lead to anisotropic Voronoi cells. Moreover, since the distribution \eqref{eq:distribution} holds for point clouds whose density does not change, the effect of changes in the local point density on the distribution of the resulting cells is uncertain.

\begin{figure}
    \centering
    \includegraphics[width=0.45\linewidth]{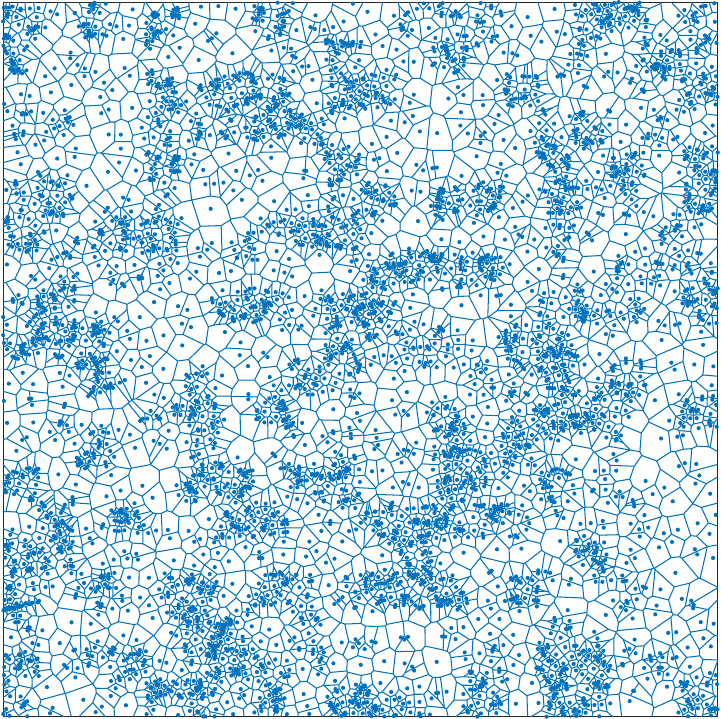} \hspace{5mm}
    \includegraphics[width=0.49\linewidth]{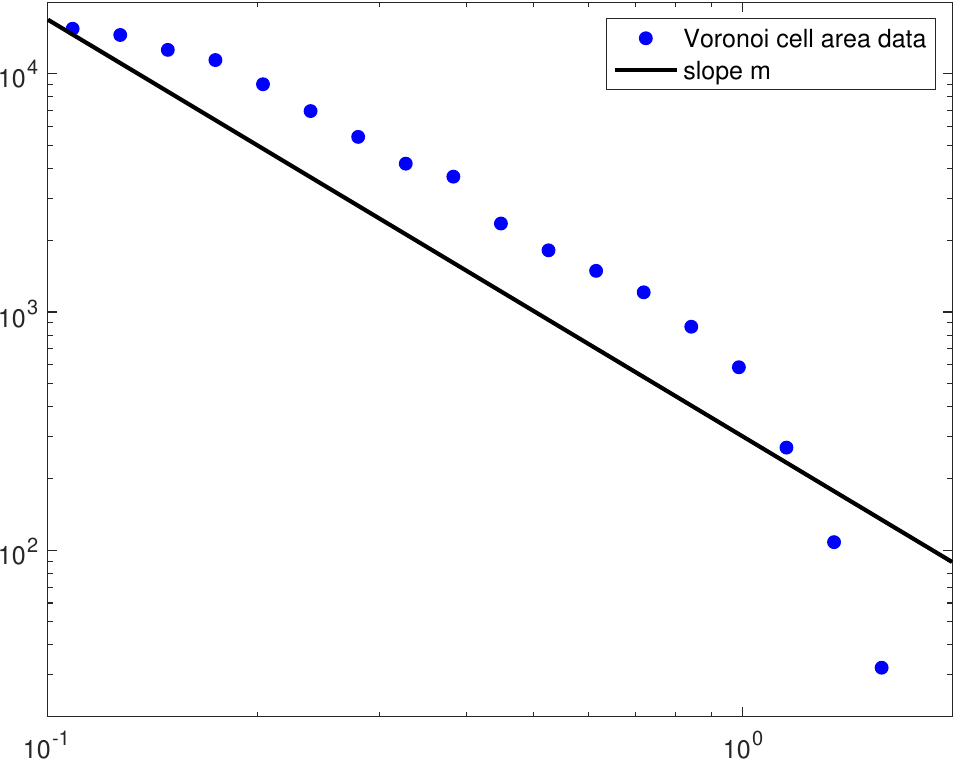}
    \caption{Voronoi cells (left) approximately satisfying the size distribution from Example 1. The empirical size distribution and the target slope are shown on the right. On the $x$-axis, we show the cell area, and on the $y$-axis the number density.}
    \label{fig:cells}
\end{figure}

\begin{figure}
    \centering
    \includegraphics[width=0.45\linewidth]{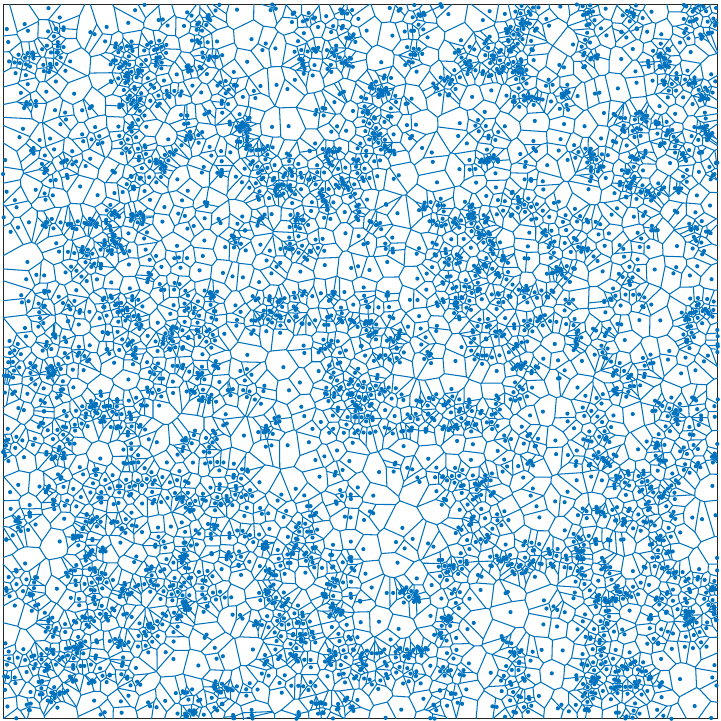} \hspace{5mm}
    \includegraphics[width=0.49\linewidth]{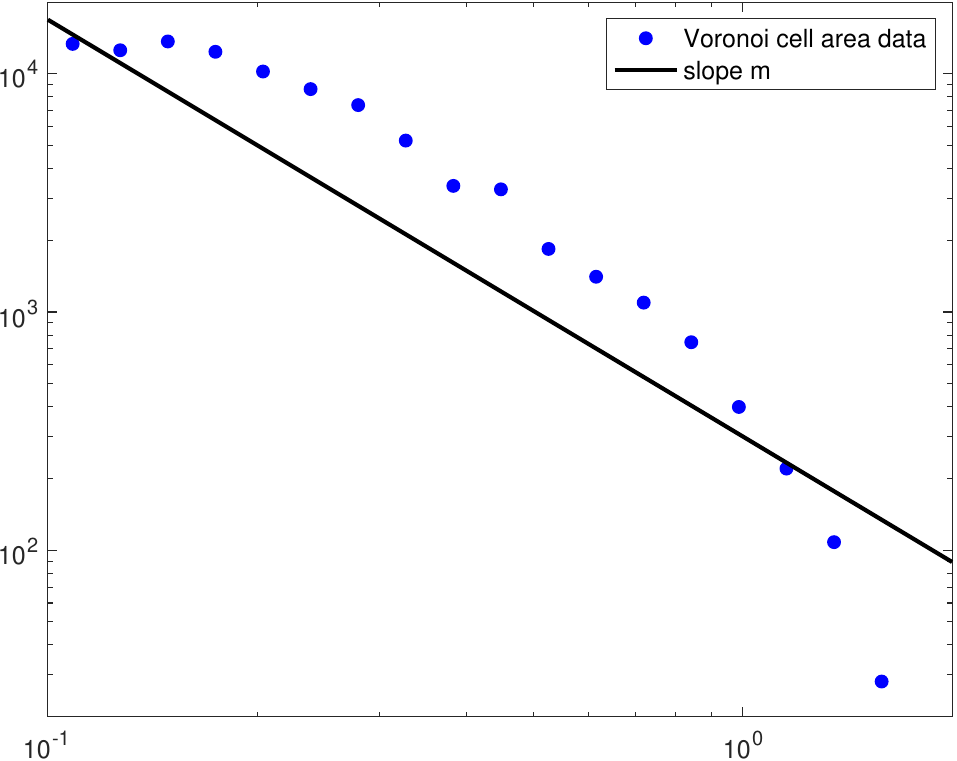}
    \caption{Same as Figure \ref{fig:cells}, but with additional mixing of the points.}
    \label{fig:cells1}
\end{figure}

We used the following practical approach to achieve a reasonable mixing of Voronoi cells of different sizes:
We divide the domain into a large number (e.g., 100 or more) of subdomains and assign to each subdomain a point density $s_i$ from the scale vector $\bs s$ such that the area fraction of subdomains with point density $s_i$ is equal to $\alpha_i$. To avoid interfaces with large point density changes, one can try to avoid adjacent subdomains corresponding to very different scales $s_i$. This approach is used to generate the cells shown on the left in Figure~\ref{fig:cells}, which corresponds to the setting of \emph{Example 1} above. Here, we split the domain into $20\times 20=400$ subdomains. Using $\bs \alpha$ from \eqref{eq:alpha_opt-ex1}, we assign $400\times 0.48$ subdomains to have point density $1/s_1$, $400\times 0.10$ subdomains to have a point density of $1/s_2$ etc. The empirical distribution closely follows the target power-law. Only the number of very large cells is too small, which is a consequence of choosing $s_1=1$ as the largest scale in $\bs s$. To improve mixing, we additionally shuffled random regions. The resulting cells and the empirical area distribution are shown in Figure \ref{fig:cells1}. We observe that the mixing has the effect that the distribution quality degrades slightly for the smallest and the largest cells. This is because the mixing makes the point cloud more uniform, thus decreasing the number of very large and very small cells.
%

\section{Conclusions}
We developed a simple method to generate Voronoi cells with a target area distribution from a point cloud. The method uses the empirical distribution of the cell areas resulting from random points with fixed density. It is based on a nonlinear least-squares problem to calculate the fractions of the domain that should correspond to random points of different density. We use a simple approach to obtain point clouds resulting in Voronoi cells that approximately satisfy the target distribution and discuss an approach to achieve a reasonable spatial mixing of cells of different sizes.

\bibliographystyle{siamplain}
\bibliography{references}

\end{document}